\documentstyle[twoside,amssymb,12pt]{article}
\def\C{{\Bbb C}}

\def\R{{\Bbb R}}
\def\P{{\Bbb P}}
\def\Z{{\Bbb Z}}

\newtheorem{prop}{Proposition}[section]
\newtheorem{dfn}[prop]{Definition}
\newtheorem{theo}[prop]{Theorem}

\newtheorem{rem}[prop]{Remark}
\newtheorem{coro}[prop]{Corollary}

\newtheorem{exam}[prop]{Example}

\begin{document}

\title{\sc Einstein-K\"ahler Metrics  on 
Symmetric Toric Fano Manifolds}  

\author{{\sc Victor V. Batyrev} \\
\small  {\em Mathematisches Institut, Universit\"at T\"ubingen}   \\
\small  {\em Auf der Morgenstelle 10,  72076  T\"ubingen, Germany}  \\
\small  {\em e-mail: batyrev@bastau.mathematik.uni-tuebingen.de} \\
and \\
{\sc Elena N. Selivanova}\thanks{Supported by Arbeitsbereich   
``Analysis'' at University of T\"ubingen}  \\
\small  {\em Nizhny Novgorod State Pedagogical University}   \\
\small  {\em Nizhny Novgorod, Russia}  \\
\small  {\em e-mail: libr@appl.sci-nnov.ru}
}
 
\date{}

\maketitle

\begin{abstract}
Let $X$ be a complex toric Fano $n$-fold and 
${\cal N}(T)$ the normalizer of a maximal torus $T$ in the 
group  of biholomorphic authomorphisms $Aut(X)$. 
We call $X$ {\em symmetric} if the trivial character is a single 
${\cal N}(T)$-invariant  algebraic character of $T$. 
Using an invariant $\alpha_G(X)$ introduced by Tian, we 
show that all  symmetric toric Fano 
$n$-folds admit an Einstein-K\"ahler metric. We remark that so far 
one doesn't know any  example of a toric Fano $n$-fold $X$ 
such that $Aut(X)$ is reductive, the  Futaki character
of $X$ vanishes, but $X$ is not symmetric. 
\end{abstract}

\thispagestyle{empty}

\newpage

\section{Introduction}

Let $X$ be a $n$-dimensional compact complex manifold with positive 
first Chern class $c_1(X)$, $g = \{  g_{i\,\overline{j}} \}$ a K\"ahler 
metric on $X$ such that the corresponding $2$-form 
\[ \omega_g = \frac{\sqrt{-1}}{2\pi} 
\sum_{i,j=1}^n g_{i\,\overline{j}} dz_i \wedge d \overline{z}_j \]
represents $c_1(X)$. It is well-known that the Ricci curvature of $g$, 
\[ Ric(g) = \frac{\sqrt{-1}}{2\pi} 
\sum_{i,j=1}^n R_{i\,\overline{j}} dz_i \wedge d \overline{z}_j, \]
\[  R_{i\,\overline{j}} = -
\frac{\partial^2 \log \det 
(g_{k\,\overline{l}})}{\partial z_i \partial  \overline{z}_j}, \]
also represents $c_1(X)$. The metric $g$ is called Einstein-K\"ahler 
if $Ric(g) =  \omega_g$.  

Let $Aut(X)$ be the  group of biholomorphic authomorphisms of $X$ and  
$Lie(Aut(X))$ the Lie algebra of $Aut(X)$. In 1957, Matsushima proved that 
if $X$ admits an Einstein-K\"ahler metric then $Aut(X)$ is a reductive 
algebraic group \cite{Mat}. 
In 1983, Futaki introduced a linear 
function $F_X\, : \,Lie(Aut(X))  \to \C$, so called {\em Futaki character}, 
which vanishes provided $X$ admits an  Einstein-K\"ahler
metric \cite{Fut}. Futaki has  conjectured  that the condition  $F_X = 0$ is 
sufficient for the existence of an Einstein-K\"ahler metric on $X$. 
Recently Tian disproved this conjecture \cite{T3}. 
This shows  that the problem of finding a sufficient condition for the 
existence of an Einstein-K\"ahler metric is rather  subtle. 

In this paper we restict ourselves to the case of 
compact complex manifolds $X$  with positive first Chern class  which 
are toric (see \cite{Dan,E,F,O}).
If $X$ is a toric Fano $n$-fold, then a maximal torus $T  \cong (\C^*)^n 
\subset Aut(X)$ has an open dense orbit $U  \cong T \subset X$. Denote 
by $M \cong \Z^n$ the group of algebraic characters of $T$. Then 
the Lie algebra $Lie(T)$ of $T$ can be identified with $N \otimes_{\Z} \C$, 
where  $N: = Hom(M, \Z)$ the dual group. Using the anticanonical 
embedding $X \hookrightarrow \P^m$ and a K\"ahler metric $g$ on $X$ 
induced by the Fubini-Study metric on $\P^m$, we obtain a 
natural moment map
\[ \mu_g\; : \; X \to M_{\R} := M \otimes_{\Z} \R \]
whose image is a convex polyhedron $\Delta$.  
The polyhedron $\Delta$ is reflexive,  and $X$ can be recovered 
from $\Delta$ as projectivization  
$X = \P_{\Delta} = Proj \, S_{\Delta}$, where 
$S_{\Delta}$ is the  graded semigroup $\C$-algebra of lattice points in the 
cone over $\Delta$ (see  \cite{Ba}). 
We denote by $R(\Delta)$ the set  of all $M$-lattice points
contained in relative interiors of codimension-$1$
faces of $\Delta$. It is well-known that 
$Aut(X)$ is reductive if and only if the set $R(\Delta)$ is centrally 
symmetric: $R(\Delta) = - R(\Delta)$. It has been shown by Mabuchi that 
if $Aut(X)$ is reductive, then 
the Futaki character $F_X$ vanishes  if and only if the barycenter 
$b(\Delta) \in M_{\R}$ 
of the polyhedron $\Delta$ is zero.  Using this  result and the complete 
classification of toric Fano $3$-folds due to the first author 
\cite{B1} and Watanabe-Watanabe \cite{WW}, Mabuchi has classified all 
Einstein-K\"ahler 
toric Fano $3$-folds. The classification  of  
$4$-dimensional  Einstein-K\"ahler toric Fano manifolds
has been obtained by Nakagawa in \cite{NY1,NY2,NY3} using results  
in \cite{B2}. Unfortunately, one toric Fano 
$4$-fold $W$  was missing in the table in \cite{B2} (see \cite{Sato}, 
Example 4.7). 
Since the Futaki character of $W$ is zero, 
one gets  a gap in the classification of Nakagawa \cite{NY2}.  
One of purposes of our paper is to fill this gap and show 
that $W$ admits an Einstein-K\"ahler metric.

Let  $X$ be a smooth projective toric  $n$-fold. Denote by 
${\cal N}(T) \subset Aut(X)$ the normalizer of a maximal torus $T$.
The group  ${\cal N}(T)$ naturally acts on $T$ by conjugations. 
This induces a linear  action of ${\cal N}(T)$ on the  group 
of algebraic characters $M = Hom_{\rm alg}(T, \C^*)$. 
Since $T$ acts trivially on $M$, the latter 
determines a linear representation of 
the finite group ${\cal W}(X):= {\cal N}(T)/T$ by integral-valued 
$n \times n$-matrices 
from $GL(M) \cong GL(n, \Z)$. 
We call $X$ {\bf symmetric}, if the trivial character 
is a single ${\cal W}(X)$-invariant (or, equivalently, 
${\cal N}(T)$-invariant)  
algebraic character of $T$: 
$$M^{{\cal W}(X)} : = \{ \chi \in M \; : \; 
\chi^g = \chi\;\;\;\mbox{\rm for all } 
g \in {\cal W}(X) \} = 0.$$  
\medskip

Our  main result is the following: 
\medskip 

\begin{theo} Let $X$ be a symmetric toric Fano $n$-fold.  
Then $X$ admits an  Einstein-K\"ahler metric.
\label{main}
\end{theo} 

It follows immediately from the definition of symmetric 
Fano manifolds that if $X = \P_{\Delta}$ is symmetric, then 
the barycenter of $\Delta$ is zero. By theorem of 
Matsushima \cite{Mat}, one also gets:

\begin{coro} 
If $X = \P_{\Delta}$ is a symmetric  toric Fano $n$-fold, then 
\[ R(\Delta) = - R(\Delta).\]
\label{roots}
\end{coro} 

It would be  interesting to know whether there exists a direct 
proof of \ref{roots} without using \ref{main}. 
We remark our theorem covers all already known examples of toric Fano 
$n$-folds ($n \leq 4$) whose  Futaki character vanish and 
whose authomorphism group is  reductive. 
It would be  interesting to know whether there exists 
an example of a toric Fano $n$-fold $X$ such that $F_X = 0$, 
$Aut(X)$ is reductive, but $X$ is not symmetric. 
Moreover, it is still unknown whether the condition $F_X = 0$ and 
$\{Aut(X)$ {\em is reductive}$\}$  is 
sufficient for the existence of an Einstein-K\"ahler
metric on toric Fano manifolds of arbitrary dimension $n$.

\bigskip

The paper is organized as follows. In Section 2 we remind 
the definition of the invariant $\alpha_G(X)$ introduced by Tian 
and its connection 
to solutions of complex Monge-Amp{\'e}re equations obtained by 
the continuity method. 
In Section 3 we give a proof of Theorem \ref{main}. 
In Section 4 we discuss several series of examples of symmetric 
toric Fano manifolds which include  all 
examples  of Einstein-K\"ahler toric Fano manifolds of dimension 
$n \leq 4$.

\bigskip
The authors would like to thank Professors Gerhard Huisken and Neil 
Trudinger for helpful discussions. The second author wish to thank 
Arbeitsbereich ``Analysis'' at University of T\"ubingen for hospitality 
and financial support. 

\section{Tian invariant $\alpha_G(X)$}

Let $X$ be a $n$-dimensional compact complex manifold with positive 
first Chern class $c_1(X)$ and  $G$ a compact subgroup of 
$Aut(X)$. Choose a $G$-invariant K\"ahler 
metric $g = \{  g_{i\,\overline{j}} \}$ on $X$ such that 
\[ \omega_g = \frac{\sqrt{-1}}{2\pi} 
\sum_{i,j=1}^n g_{i\,\overline{j}} dz_i \wedge d \overline{z}_j \]
represents $c_1(X)$. One has a natural $G$-invariant 
volume form $dV_g$ on $X$
$$dV_g := \frac{\omega_g^n}{n!}, \;\;\;\; Vol_g(X) := \int_X dV_g = 
\frac{c_1^n(X)}{n!}.$$  

It is well-known that the problem of finding an Einstein-K\"ahler 
metric on $X$ is equivalent to solving the following complex 
Monge-Amp{\'e}re equation for smooth real-valued functions $\varphi$ on $X$:
\begin{equation}
\det \left(  g_{i\,\overline{j}} + 
\frac{\partial^2 \varphi }{\partial z_i \partial  \overline{z}_j} \right)
=  \det (g_{i\,\overline{j}})  e^{F - t\varphi}, \;\;\; \forall t \in [0,1]
\label{m-a} 
\end{equation}
where the smooth real-valued function $F$ is defined 
by the conditions:
\[ \frac{\partial^2 \varphi }{\partial z_i \partial  \overline{z}_j} = 
 R_{i\,\overline{j}} -  g_{i\,\overline{j}}, \;\; 
\int_X e^F dV_g = Vol_g(X). \]
If $\varphi$ is a solution of  (\ref{m-a}) for $t =1$, then 
\[  g_{i\,\overline{j}}' :=    g_{i\,\overline{j}} + 
\frac{\partial^2 \varphi }{\partial z_i \partial  \overline{z}_j} \]
is an Einstein-K\"ahler metric on $X$. By famous  theorem of Yau,
there exists always a solution of  (\ref{m-a}) for all  
$t \in [0, \varepsilon)$ if  $\varepsilon$ is sufficiently small. 
Using the continuity method, one can show that the existence 
of a solution $\varphi$ for $t =1$ is equivalent to zero-order  
{\em a priori} estimates of $\varphi$. 

Let us recall the definition of an invariant $\alpha_G(X)$ introduced 
by Tian \cite{T1}:

\begin{dfn} 
{\rm Let  $P_G(X, g)$ be the set of all $C^{2}$-smooth $G$-invariant 
real-valued functions $\phi$ such that $\sup_X \phi = 0$ and 
\[ \omega_g  +   \frac{\sqrt{-1}}{2\pi} \partial \overline{\partial} \phi \]
is a nonnegative $(1,1)$-form.
Then {\bf Tian invariant} ${\alpha}_G(X)$ 
is defines as supremum of all $\lambda >0$ such that 
\[ \int_X e^{- \lambda \phi} dV_g \leq C(\lambda)\;\;\; \forall 
\phi \in P_G(X, g),  \]
where $C(\lambda)$ is a positive constant depending only on $\lambda$, 
$g$ and $X$. } 
\end{dfn} 

\begin{rem}
{\rm 
It is easy to show that $\alpha(X)$ doesn't depend on the choice of a 
$G$-invariant metric $g$. Moreover, $\alpha_G(X)$ doesn't change if in the 
above definition we replace $P_G(X, g)$ by a smaller subset consisting 
of all $C^{\infty}$-smooth $G$-invariant real-valued functions $\phi$ 
such that $\sup_X \phi = 0$ and 
\[ \omega_g  +   \frac{\sqrt{-1}}{2\pi} \partial \overline{\partial} \phi \]
is a positive definite $(1,1)$-form (see \cite{Tian}). } 
\end{rem} 

Deriving a zero-order 
{\em a priori} estimate for the solutions of $( \ref{m-a})$, 
Tian has proved the following important result (\cite{T1},  
Theorems 2.1 and 4.1):

\begin{theo}
Let $X$ be a Fano $n$-fold and $G \subset Aut(X)$ is a compact 
subgroup such that $$\alpha_G(X) > \frac{n}{n+1}.$$
Then $X$ admits an Einstein-K\"ahler metric. 
\label{invar}
\end{theo}

\section{Main theorem}

Throughout this section we use standard notations from  the theory of 
toric varieties (see e.g. \cite{Dan}). Let $M$ be a free abelian group of 
rank $n$, $N = Hom(M,\Z)$ the dual group, $M_{\R}:= M \otimes_{\Z} {\R}$ 
 $N_{\R}:= N \otimes_{\Z} {\R}$. Denote by 
$\langle *, * \rangle\, : \, M_{\R} \times N_{\R} \to \R$ the canonical 
nondegenerate pairing.  
Let $X = X_{\Sigma}$ be a smooth projective toric $n$-fold defined by 
a complete fan $\Sigma$ of regular cones $\sigma \subset N_{\R}$. Then 
a maximal torus $T \subset Aut(X)$ acting  on $X$ has an open dense orbit 
$U \subset X$.  The normalizer ${\cal N}(T) \subset Aut(X)$ of $T$ 
has a natural action on $U$. 
Let us set  ${\cal W}(X):= {\cal N}(T)/T$. 
By functorial properties of toric varieties (see 
\cite{Dan}, \S 5), one 
immediately obtains:

\begin{prop} 
Let $X= X_{\Sigma}$ be a smooth projective toric 
$n$-fold defined by a complete
regular  polyhedral fan $\Sigma$.  Then 
the group ${\cal W}(X)$ is isomorphic to the finite group of all symmetries 
of $\Sigma$, i.e.,  ${\cal W}(X)$ is isomorphic to a subgroup 
of $GL(M)$ $(\cong GL(n, \Z))$ consisting of all elements 
$\gamma \in GL(M)$ such that $\gamma(\Sigma) = \Sigma$. 
\label{group}
\end{prop} 

Since the open subvariety 
$U \subset X$ is a principal homogeneous space of $T$, 
we can identify   $U$ with $T$ by choosing an arbitrary point 
$x_0 \in U$. This indentification defines a splitting 
of the short exact sequence
\[ 1 \to T \to {\cal N}(T) \to {\cal W}(X) \to 1, \]   
i.e., an embedding ${\cal W}(X) \hookrightarrow {\cal N}(T) \subset Aut(X)$. 
We denote by ${\cal W}(X,x_0)$ the image of ${\cal W}(X)$ in $Aut(X)$  
under this embedding. 
Denote by ${\cal K}(T) \cong (S^1)^n$ the maximal
compact subgroup in $T$. In the sequel we shall use the canonical 
isomorphism $T / {\cal K}(T) \cong N_{\R}$ and the isomorphism 
$U / {\cal K}(T) \cong   N_{\R}$ which identifies the orbit 
${\cal K}(T) x_0$ with the zero element $0 \in N_{\R}$. The last isomorphism
shows that the  ${\cal N}(T)$-action on $U$ descends to a linear action 
of  ${\cal W}(X)$ on  $N_{\R}$. If one chooses an integral  basis
$e_1, \ldots, e_n$ of $N$ and the dual basis 
$e_1^*, \ldots, e_n^*$ of $M$, then the  induced  isomorphisms 
$N_{\R} \cong \R^n$, $M_{\R} \cong \R^n$ and $T \cong (\C^*)^n$ 
allow to introduce affine logarithmic coordinates $y_i = \log |z_i|$ $(i =1, 
\ldots, n)$ on $N_{\R}$,  where $z_1, \ldots, z_n$ the 
standard holomorphic coordinate system on $(\C^*)^n$.
We choose  $G$ to be the maximal compact subgroup 
in ${\cal N}(T)$  generated by ${\cal W}(X,x_0)$ and  ${\cal K}(T)$, so that  
we have   the short exact sequence
\[  1 \to {\cal K}(T) \to G \to {\cal W}(X) \to 1. \]

Now we assume that a projective toric $n$-fold $X$ has positive first
Chern class. In this case,  one  obtains a convex  ${\cal W}(X)$-
invariant polyhedron $\Delta \subset M_{\R}$ defined by the affine linear 
inequalities $\langle y, e \rangle \leq 1$ where $e$ runs over 
all primitive integral 
generators $e$ of $1$-dimensional cones $ \sigma = \R_{\geq 0}e 
\in \Sigma$. 
Let $L(\Delta) = \{ v_0, v_1, \ldots, v_m \} := M \cap \Delta$. 
Then   
$v_0, v_1,  \ldots, v_l$ determine  algebraic characters 
$\chi_i\, : \,  T \to \C^*$ of $T$ $(i =0, \ldots, m)$.  
Moreover, we have 
\[ |\chi_i(x) | = e^{\langle v_i, y \rangle }, \;\; i =0, \ldots, m, \]
where $y$ is the image of $x$ under the canonical  projection 
$\pi\, : \, T \to N_{\R}$.  
Let us define the function $u\, : \, U \to \R$ as follows:
\begin{equation} u :=  
\log( \sum_{i =0}^m |\chi_i(x)| ), \;\; x \in U \cong T. 
\label{u-def} 
\end{equation} 
Since $u$ is  ${\cal K}(T)$-invariant, $u$ descends to a 
function $\tilde{u}\, : \, N_{\R} \to \R$ defined as
\begin{equation} \tilde{u} :=  
\log( \sum_{i =0}^m e^{\langle v_i, y \rangle} ), \;\; y \in N_{\R}. 
\label{u-def2} 
\end{equation}  
Since 
$L(\Delta)$ is ${\cal W}(X)$-invariant, 
one obtains  the following $G$-equivariant  moment map
\[ \mu_{\tilde{u}}\; : \; N_{\R} \to M_{\R} \]
$$ y = (y_1, \ldots, y_n) 
\mapsto {\rm Grad}\, \tilde{u} 
:= \left(\frac{\partial \tilde{u} }{\partial y_1}(y),
 \cdots, 
\frac{\partial \tilde{u}}{\partial y_n}(y) \right)$$
which is a diffeomorphism of $N_{\R}$ with the interior 
of the polyhedron $\Delta$. 

Consider the  $G$-invariant hermitian metric $g =     
\{  g_{i\,\overline{j}} \}$ on $X$ such that the restriction 
of the  corresponding to $g$  
differential $2$-form on $U$ is defined  by 
$$\omega_g = \frac{\sqrt{-1}}{2\pi} 
\partial \overline{\partial} u.$$
We remark that  the metric  $g$ is exactly the pull-back 
of the Fubuni-Study metric form $\P^m$ with respect to the anticanonical 
embedding $X \hookrightarrow \P^m$ defined by the algebraic characters 
$\chi_0, \chi_1, \ldots, \chi_m$. Then the restriction of the moment 
$\mu_g\,: \, X \to M_{\R}$ to $U$ is exactly the composition 
of the canonical projection 
$\pi\, : \, T \to N_{\R}$ and $\mu_{\tilde{u}}\, : \, N_{\R} \to 
M_{\R}$. In particular, $\Delta = \mu_g(X)$.

Using the above considerations, one can derive from the complex  
Monge-Amp{\'e}re equation (\ref{m-a}) for 
a $G$-invariant function $\varphi\, : \, X \to \R$ 
the  real Monge-Amp{\'e}re equation
\begin{equation} 
 \det \left( \frac{ \partial^2 (\tilde{u} + \tilde{\varphi})}{\partial y_i 
\partial y_j} 
\right)  = \exp ( -\tilde{u} - t \tilde{\varphi}),\; \; \forall t \in 
[0,1],
\label{m-a-r}
\end{equation}
where $\tilde{\varphi}$ is a smooth ${\cal W}(X)$-invariant 
real-valued function on $N_{\R}$ obtained   as descent of $\varphi |_U$ 
to $N_{\R}$.

\begin{prop} 
Let $X$ be a toric Fano $n$-fold with $G$-action as above. Denote by  
$dy$  the volume $n$-form 
on $N_{\R}( \cong \R^n)$ corresponding to the Haar measure 
on $N_{\R}$ normalized by the lattice 
$N \subset N_{\R}$. 
Let  $\tilde{\alpha}_G(X)$ be 
the supremum of all $\lambda >0$ such that 
\[ \int_{N_{\R}} e^{- \lambda \tilde{\phi} -\tilde{u} } dy 
\leq \tilde{C}(\lambda)\;\;\; \forall 
\tilde{\phi} \in P_G(N_{\R},\tilde{u}), \]
where $P_G(N_{\R}, \tilde{u} )$ is  the set of all $C^{2}$-smooth 
${\cal W}(X)$-invariant 
functions $\tilde{\phi}\, : \, N_{\R} \to \R$ such that $\tilde{u} + 
\tilde{\phi}$ is 
upper convex,  $\sup_X \tilde{\phi} = 0$, and 
$|\tilde{\phi}|$ is bounded on the whole $N_{\R}$. 
Then 
\[\tilde{\alpha}_G(X)  \leq {\alpha}_G(X). \]
\label{comp}
\end{prop}

\noindent
{\em Proof.} 
 Let $\phi$ be an element of $P_G(X,g)$. Since $\phi$ is 
${\cal K}(T)$-invariant, the restriction of $\phi$ to $U$ 
 descends to a smooth 
$C^2$-function real-valued $\tilde{\phi}$ on 
$N_{\R} \cong U/{\cal K}(T)$. Moreover, 
it follows from  $G$-variance of $\phi$ that  $\tilde{\phi}$ is 
invariant under the finite group ${\cal W}(X)$ 
acting linearly on $N_{\R}$. 
The nonnegativity of  the $(1,1)$-form  
\[ \omega_g  +   \frac{\sqrt{-1}}{2\pi} \partial \overline{\partial} 
{\phi} = 
  \frac{\sqrt{-1}}{2\pi} \partial \overline{\partial} ( u + \phi) \]
immediately implies that the matrix 
\[ \left( \frac{ \partial^2 (\tilde{u} + 
\tilde{\phi} )}{\partial y_i \partial y_j} 
\right) 
\]
is nonnegative definite, i.e., $\tilde{u} + \tilde{\phi}$ 
is an upper convex function on $N_{\R}$. 
Let $d\theta$ be  a volume $n$-form defining the 
canonically normalized Haar measure on the compact group ${\cal K}(T)$.
We remark that the restriction of the 
volume $2n$-form  $dV_g$ to $U \cong T$ 
equals   $h e^{-u} dy d\theta$, where $h$ is a smooth real-valued bounded 
function on $X$. Therefore, the inequality 
\[ \int_X e^{- \lambda \phi} dV_g \leq C(\lambda)\;\;\; \forall 
\phi \in P_G(X, g)  \]
immediately follows from 
\[ \int_{N_{\R}} e^{- \lambda \tilde{\phi} - \tilde{u}} dy 
\leq \tilde{C}(\lambda)\;\;\; \forall 
\tilde{\phi} \in P_G(N_{\R}, \tilde{u}). \]
Thus, we have $\tilde{\alpha}_G(X)  \leq {\alpha}_G(X)$.
\hfill $\Box$

\begin{prop} Let $X = \P_{\Delta}$ be a toric Fano $n$-fold and 
$\tilde{u}$ the  function defined by (\ref{u-def2}). Choose 
$\tau$ to be an arbirary positive real number. Then 
\[  \int_{N_{\R}} e^{- \tau \tilde{u}} dy   \leq \frac{ v(\Delta)}{\tau^n}, \]
where $v(\Delta)$ is the number 
of vertices of $\Delta$.   
\label{e-int}
\end{prop} 

\noindent
{\em Proof.} Let $v(\Delta) = l$. Denote by  $w_1, \ldots, w_l$ 
all vertices of $\Delta$. It follows from the formula (\ref{u-def2})
that for all $y \in N_{\R}$ we have 
\begin{equation}
\tilde{u}(y) > \langle w_j, y \rangle, \;\; j =1, \ldots, l,  
\end{equation}  
and hence $\tilde{u}(y) > \overline{u}(y)$, 
where 
$\overline{u}  := \max_{j =1, \ldots,l}  \langle w_j, y \rangle$. 
Therefore, we obtain
\begin{equation}
\int_{N_{\R}} e^{- \tau \tilde{u}} dy \leq   
\int_{N_{\R}} e^{- \tau \overline{u}} dy. 
\label{estim}
\end{equation}
It follows from definition of $\Delta$ 
that $l$ is exactly the number of $n$-dimensional cones $\sigma_1, \ldots, 
\sigma_l$ in the fan $\Sigma$ defining $X$. Moreover,  $\overline{u}$ 
is a continuous piecewise linear function whose restriction 
to $\sigma_j$ equals $\langle w_j, y \rangle$.  
On the other hand, 
\begin{equation} 
\int_{\sigma_j} e^{- \tau \overline{u}} dy = 
\int_{\R_{\geq 0}^n} e^{- \tau (y_1 + \cdots + y_n)} dy_1 \cdots dy_n = 
\prod_{i =1}^n \left( \int_{\R_{\geq 0}} e^{- \tau y_i } dy_i  \right) = 
\frac{1}{\tau^n},
\label{1-int}   
\end{equation} 
since every $n$-dimensional cone $\sigma_j \in \Sigma$ $(j =1, \ldots, l)$ 
 is generated by a basis of the lattice $N$. 
Using $N_{\R} = \sigma_1 \cup \cdots \cup \sigma_l$ together 
with  (\ref{estim}) and  (\ref{1-int}), we come to the required inequality. 
\hfill $\Box$
\bigskip

The next statement plays the crucial role in the proof of Theorem \ref{main}: 

\begin{theo} 
Let $X= \P_{\Delta}$ be a symmetric toric Fano $n$-fold and $\tilde{\phi}$ 
is an arbitrary function from $P_G(N_{\R}, \tilde{u})$. Then 
\[  \tilde{u}(y)   + \tilde{\phi}(y) \geq 0 \;\;\;\; \forall y \in N_{\R}. \]
\label{pos}
\end{theo} 

\noindent
{\em Proof.} Let $\tilde{\phi}$ be  
an arbitrary function from $P_G(N_{\R}, \tilde{u})$. Consider the following 
moment map: 
\[ \mu_{\tilde{u} + \tilde{\phi}}\; : \; N_{\R} \to M_{\R}, \]
$$ y = (y_1, \ldots, y_n) 
\mapsto {\rm Grad}\,(\tilde{u} + \tilde{\phi})(y)
 := \left(\frac{\partial (\tilde{u} +  \tilde{\phi})
  }{\partial y_1}(y),
 \cdots, 
\frac{\partial (\tilde{u} + \tilde{\phi})  }{\partial y_n}(y) \right).$$

First of all we show that $\mu_{\tilde{u} + \tilde{\phi}}(N_{\R}) \subset 
\Delta$. Let $z = \mu_{\tilde{u} + \tilde{\phi}}(y')$ for 
some $y' \in N_{\R}$. It follows from the convexity of 
$\tilde{u} + \tilde{\phi}$  
that  for all $y \in N_{\R}$ one has 
$$\tilde{u}(y) + \tilde{\phi}(y) \geq \langle z, y - y' \rangle + 
\tilde{u}(y') + \tilde{\phi}(y').$$
In other words, the function 
 $\tilde{u}(y) + \tilde{\phi}(y) - \langle z, y \rangle$
attains the  global minimum at  $y' \in N_{\R}$.      
Let $\{ w_1, \ldots, w_l \}$ be the set of all vertices of $\Delta$ and 
$\overline{u}:=\max_{j =1, \ldots,l}  \langle w_j, y \rangle$ the 
piecewise linear function as in the proof of \ref{e-int}.    
Using obvious inequalities
\[ \log l  + \overline{u} \geq \tilde{u} \geq  \overline{u} \]
and the fact that  $\tilde{\phi}$ is globally bounded on $N_{\R}$, 
we conclude that the piecewise linear function 
$\overline{u}(y) -  \langle z, y \rangle$ is 
bounded from below on the whole $N_{\R}$. The latter is possible only 
if  $\overline{u}(y) -  \langle z, y \rangle \geq 0$ for all $y \in N_{\R}$. 
Since $\overline{u}(e) =1$ for all primitive integral generators of 
$1$-dimensional cones $\sigma \in \Sigma$, we obtain that 
for all these generators holds  $\langle z, e \rangle \leq 1$, i.e., 
$z \in \Delta$. 

Since $\sup_{N_{\R}} \tilde{\phi} =0$, there exists a sequence 
$\{ q_k \}_{k \geq 1} $ of points $q_k \in N_{\R}$ 
such $-1/k \leq  \tilde{\phi}(q_k) \leq 0$. Denote $z_k =  
\mu_{\tilde{u} + \tilde{\phi}}(q_k)$. Since all $z_k$ belong to $\Delta$, 
we can assume without loss of generality that 
\[ \lim_{k \to \infty} z_k = z \in \Delta \] 
(otherwise one chooses an appropriate subseqeunce of $\{ q_k \}_{k \geq 1}$).
It follows from the convexity of  
$\tilde{u} + \tilde{\phi}$  
that  for all $y \in N_{\R}$ and all $k \geq 1$ one has 
$$\tilde{u}(y) + \tilde{\phi}(y) - \langle z_k, y \rangle  \geq 
\tilde{u}(q_k) + \tilde{\phi}(q_k) - \langle z_k, q_k \rangle .$$

Now we remark that $\tilde{u}(q_k) \geq \overline{u}(q_k)  \geq  
\langle z_k, q_k \rangle $ for 
all $k \geq 1$, because   $z_k$ is contained in $\Delta$. Therefore, we have 
\[ \tilde{u}(y) + \tilde{\phi}(y) - \langle z_k, y \rangle  \geq -1/k, 
 \;\; \; \forall
y \in N_{\R} .\]
Taking limit $k \to \infty$, we obtain
\begin{equation}
\tilde{u}(y) + \tilde{\phi}(y) - \langle z, y \rangle  \geq 0, \;\; \; \forall
y \in N_{\R}.
\label{geq0}
\end{equation}

We set  $r := |{\cal W}(X)|$ and consider the points  
$z, \gamma_1 z, \ldots, \gamma_{r-1}z$, where $\{ 
\gamma_1, \ldots, \gamma_{r-1} \}$
the set of all  elements of ${\cal W}(X)$ which are different from the 
identity. Since $\tilde{u} + \tilde{\phi}$  is ${\cal W}(X)$-invarint,
we obtain from (\ref{geq0}) $r-1$ additional  inequalities:
\begin{equation}
\tilde{u}(y) + \tilde{\phi}(y) - \langle \gamma_j z, y \rangle  \geq 0, 
\;\;\; \;\; \; \forall
y \in N_{\R}, \;\; j =1, \ldots, r-1.
\label{geq0-add}
\end{equation}  
Now we remark that 
\[ z' := z + \sum_{j =1}^{r-1} \gamma_j z \]
is obviously ${\cal W}(X)$-invariant. Using the fact that $X$ is 
a symmetric Fano $n$-fold, we conclude that $z' = 0$. Summing the inequalities 
in (\ref{geq0}) and (\ref{geq0-add}), we obtain 
\[  \tilde{u}(y)   + \tilde{\phi}(y) \geq 0 \;\;\;\; \forall y \in N_{\R}. \]
\hfill $\Box$
\bigskip

\noindent
{\bf  Proof of Theorem \ref{main}.} Choose arbitrary $\lambda \in (0,1)$ and 
$\tilde{\phi} \in P_G(N_{\R}, \tilde{u})$. Using \ref{pos} and \ref{e-int}, 
we obtain  
\[ \int_{N_{\R}} e^{-\lambda \tilde{\phi} -\tilde{u}}dy = 
\int_{N_{\R}}  e^{-\lambda  (\tilde{\phi} + \tilde{u}) }  
e^{(\lambda-1) \tilde{u} }dy \leq
\sup_{N_{\R}} \left\{  e^{-\lambda (\tilde{\phi} + \tilde{u})} \right\} 
\int_{N_{\R}} e^{(\lambda -1) \tilde{u} }dy 
\leq \frac{v(\Delta)}{(1- \lambda)^n}.  \]
Therefore, $\tilde{\alpha}_G \geq 1$. By  
\ref{comp} and \ref{invar}, we conclude that 
$X$ admits an Einstein-K\"ahler metric.  
\hfill $\Box$ 
\bigskip

\section{Some examples} 

In this section we consider series of examples of   
symmetric toric Fano $n$-folds which include many already known 
examples of toric Einstein-K\"ahler manifolds. 

\begin{exam} 
{\rm Let $V_{k}$ smooth projective toric Fano $n$-fold $(n = 2k)$ defined 
by a  fan $\Sigma$ of regular polyhedral cones whose generators are 
$\pm e_1, \ldots, \pm e_n, \pm ( e_1 + \cdots + e_n)$, 
where $e_1, \ldots, e_n$ is an integral basis of the lattice $N$. 
The toric Fano $n$-fold $V_k$  has been introduced 
by Voskresensky and Klyachko \cite{VK}. 
Since the corresponding polyhedron $\Delta = 
\Delta(V_{k})$ is centrally symmetric, 
$V_{k}$ is a symmetric toric Fano $n$-fold (see \ref{group}). We remark that  
$V_{1}$ is $\P^2$ with $3$ points blown-up. The existence of an 
Einstein-K\"ahler metric on $V_1$ was proved by Siu \cite{S1}, Tian-Yau 
\cite{TY}, and Nadel \cite{N}.  The existence of 
 an Einstein-K\"ahler metric on the $4$-fold $V_{2}$ was proved 
by Nakagawa in \cite{NY1} using results of Nadel \cite{N}. 
}
\end{exam}

\begin{exam}
{\rm Let $k,m$ be integers satisfying the condition 
$1 \leq k \leq m$. Denote by $S_{m,k}$ toric Fano $n$-fold 
$(n = 2m +1)$ which is the projectivization $\P(E)$ of the 
split bundle   $E = {\cal O} \oplus {\cal O}(k, -k)$ over 
$\P^m \times \P^m$.  This toric manifold is defined 
by a fan $\Sigma$ whose cones have  the following $2m + 4$ generators: 
\[ e_1, \ldots, e_{2m}, \pm e_{2m+1},  \; 
-(e_1 + e_2 + \cdots + e_m + ke_{2m+1}), \] \[ 
 -(e_{m+1} + e_{m+2} + 
\cdots + e_{2m} - ke_{2m+1}),  \]
where  $e_1, \ldots, e_{2m+1}$ is an integral basis of $N$. 
There exist an   authomorphisms $\alpha$ of $\Sigma$ of order $m+1$ 
such that 
\[ \alpha(e_{2m+1})=e_{2m+1}, \; \; 
\alpha(e_i) = e_{i+1}, \; \; \alpha(e_{i+m}) = e_{i+m+1},  \;\; 
i =1, \ldots, m-1;  \]  
\[ \alpha(e_m) = - (e_1 + \ldots + e_m  + ke_{2m+1}), \;\; 
 \alpha'(e_{2m}) = - (e_{m+1} + \ldots + e_{2m}- ke_{2m+1}). \]
There exists an authomorphism $\beta$ of order $2$ defined by 
\[ \beta(e_{2m+1}) =  -e_{2m+1},\;   \beta(e_i) = e_{i +m}, 
\; \beta(e_{i+m}) = e_{i} \; 
\;\;i =1, \ldots, m. \]  
The common fix point set of $\alpha$ and $\beta$ is exactly  
$0 \in N_{\R}$. By \ref{group},   $S_{m,k}$ is a symmetric 
toric Fano $n$-fold.

The Einstein-K\"ahler manifold $S_{m,k}$ was discovered 
by Sakane \cite{S}. The existence 
of an Einstein-K\"ahler metric on $S_{m,k}$ was obtained by Mabuchi  
using another method (see (10.3.2) in  \cite{M}).  
We  remark that $S_{m,1}$ is isomorphic to $\P^{2m +1}$ blown-up at two skew 
$m$-dimensional subspaces. The existence of an 
Einstein-K\"ahler metric on  $S_{m,1}$ was proved independently  
by Nadel (\cite{N}, Example 6.4).} 
\end{exam}

\begin{exam} 
{\rm Choose integers $k, m$ such that 
$0 \leq k \leq m$.  In \cite{NY2} Nakagawa introduced 
a toric Fano $n$-fold $X_{m,k}$  
$(n = 2m + 2)$ defined by a fan $\Sigma$ whose $2m +8$ generators are 
\[  e_1,  \ldots, e_{2m}, \pm e_{2m+1}, \pm e_{2m+2}, 
\pm (e_{2m+1} + e_{2m+2}), \]
\[ -(e_1 + \cdots + e_m  - ke_{2m+1}),  \;\; 
-(e_{m+1} + \cdots + e_{2m}  + ke_{2m+1}). \]
There exist an authomorphism $\alpha$ of $\Sigma$ of order $m+1$ 
such that 
\[ \alpha(e_{2m+1}) = e_{2m+1}, \;\; \alpha(e_{2m+2}) = e_{2m+2}, \;\; \]
\[ 
\alpha(e_i) = e_{i+1}, \;\;  \alpha(e_{i+m}) = e_{i+m+ 1}, \;\; 
i =1, \ldots, m-1;  \]
\[  \alpha(e_m) = - (e_1 + \ldots + e_m -k e_{2m+1}), \;\; 
\alpha(e_{2m}) = - (e_{m+1} + \ldots + e_{2m} + ke_{2m+1}). \]
On the other hand, there exists  an authomorphism 
$\beta$ of $\Sigma$ of order $2$ 
defined by  
\[ \beta(e_i) = e_{i+m}, \;\;  \beta(e_{i+m}) = e_{i}, \;\; 
i =1, \ldots, m; \]
\[ \beta(e_{2m+1}) = - e_{2m+1}, \;\; \beta(e_{2m+2}) = -e_{2m+2}. \]
The common fix point set of $\alpha$ and $\beta$ is exactly  
$0 \in N_{\R}$. By \ref{group},  $X_{m,k}$ is a symmetric toric Fano 
$n$-fold.  
 The existence of 
 an Einstein-K\"ahler metric on $X_{m,k}$  was proved 
by Nakagawa in \cite{NY1} using results of Nadel \cite{NY2}. 

} 
\end{exam} 

\begin{exam} 
{\rm Let $W_{m}$ be  $\P^m \times \P^m$ blown-up 
along $m +1$ codimension-$2$ subvarieties $Z_i \cong \P^{m-1} 
\times \P^{m-1}$  defined by the equations 
$z_i =0$, $z_i' = 0$ $(i =0,1, \ldots, m)$, 
where $(z_0: z_1: \cdots : z_m)$ and $(z_0':  z_1': \cdots : z_m')$ are 
homogeneous coordinates on two $\P^m$'s. The toric manifold 
$W_m$ is determined by a $2m$-dimensional  
fan $\Sigma \subset N_{\R}$ whose cones have the following $3m + 3$ generators 
\[ e_1, \ldots, e_{2m}, -(e_1 + \ldots + e_m),  
-(e_{m+1} + \ldots + e_{2m}), -(e_1 + \ldots + e_{2m}), \]
\[ e_i + e_{i+m}, \;\; i =1, \ldots, m, \]
where $e_1, \ldots, e_{2m}$ is an integral basis of $N$.    
There exists an authomorphism $\alpha$ of $\Sigma$ of order $m+1$ 
such that
\[ \alpha(e_i) = e_{i+1}, \;\;  \alpha(e_{i+m}) = e_{i+m+ 1}, \;\; 
i =1, \ldots, m-1,  \]
\[  \alpha(e_m) = - (e_1 + \ldots + e_m), \;\; 
\alpha(e_{2m}) = - (e_{m+1} + \ldots + e_{2m}). \]
On the other hand, there exists  an authomorphism 
$\beta$ of $\Sigma$ of order $2$ 
defined by  
\[ \beta(e_i) = e_{i+m}, \;\;  \beta(e_{i+m}) = e_{i}, \;\; 
i =1, \ldots, m. \]
The common fix point set of $\alpha$ and $\beta$ is exactly  
$0 \in N_{\R}$.  By \ref{group}, $W_m$ 
is a symmetric toric Fano $n$-fold $(n = 2m)$.

We remark that $W_1 = V_1$ is again  $\P^2$ with $3$ points blown-up.
The toric Fano $4$-fold $W_2$ is exactly the single one
missed in the table 
\cite{B2}. In particular, we come to conclusion that there exist 
exactly $12$ different Einstein-K\"ahler toric Fano $4$-folds 
(cf. \cite{NY2,NY3} and \cite{Sato}, Example 4.7).   
} 
\end{exam}

We remark that  any  Einstein-K\"ahler toric Fano manifold $X$  of 
dimension $n \leq 4$ which
can not be decomposed into a product of lower dimensional varieties 
is  either a projective space, or one of the toric Fano manifolds from 
4.1-4.4.


\begin{thebibliography}{99}

\bibitem{B1} V. Batyrev, {\em 
Toroidal Fano 3-folds}, Math. USSR, Izv. {\bf 19} (1982), 13-25.

\bibitem{Ba}  V. Batyrev, {\em Dual polyhedra and mirror symmetry
 for Calabi-Yau hypersurfaces in toric varieties}, 
J. Alg. Geom. {\bf 3}, No.3 (1994), 493-535.

\bibitem{B2}  V. Batyrev, {\em  On the classification of toric 
Fano $4$-folds}, math.AG/9801107. 

\bibitem{Dan} V. I. Danilov, {\em Geometry of toric varieties}, 
Russ. Math. Surv. {\bf 33}, No.2 (1978), 97-154.  

\bibitem{E} G. Ewald, {\em Combinatorial convexity and algebraic geometry}, 
Graduate Texts in Mathematics, {\bf 168}, New York, Springer (1996).

\bibitem{F} W. Fulton, {\em Introduction to Toric Varieties}, 
Ann. of Math. Studies {\bf 131}, Princeton Univ. Press, 1993. 

\bibitem{Fut} A. Futaki, {\em An obstruction to the existence 
of Einstein K\"ahler metrics}, 
Invent. Math. {\bf 73} (1983), 437-443.

\bibitem{M} T. Mabuchi, {\em Einstein-K\"ahler forms, Futaki 
invariants and convex geometry on toric Fano varieties}, 
Osaka J. Math. {\bf 24} (1987), 705-737.

\bibitem{Mat} Y. Matsushima,  {\em 
Sur la structure du groupe d'hom{\'e}omorphismes analytiques d'une 
certaine vari{\'e}t{\'e} kaehl{\'e}rienne},  Nagoya Math. J. {\bf 11} 
(1957), 145-150. 

\bibitem{N} A. Nadel, {\em Multiplier ideal sheaves and K\"ahler-Einstein 
metrics of positive scalar curvature},  
Ann. Math., II. Ser. 132, No.3  (1990), 549-596.  

\bibitem{NY1} Y. Nakagawa, {\em  
Einstein-K\"ahler toric Fano fourfolds}, 
Tohoku Math. J., II. Ser. 45, No.2 (1993), 297-310. 

\bibitem{NY2} Y. Nakagawa, {\em 
Classification of Einstein-K\"ahler toric Fano fourfolds},  
Tohoku Math. J., II. Ser. 46, No.1 (1994),  125-133.

\bibitem{NY3} Y. Nakagawa, {\em Combinatorial Formulae for 
Futaki Characters and Generalized Killing Forms on Toric Fano
Orbifolds}, Preprint 1997. 

\bibitem{O} T. Oda, {\em Convex Bodies and Algebraic Geometry - An 
introduction to the thoery of toric varieties}, Ergebnisse Math. Grenzgeb. 
(3), Vol. 15, Springer-Verlag, Berlin, Heidelberg, New York, London, Paris,
Tokyo, 1988. 

\bibitem{S} Y. Sakane, {\em 
Examples of compact Einstein K\"ahler manifolds with positive Ricci 
tensor}, Osaka J. Math. {\bf 23} (1986), 585-616. 

\bibitem{Sato} H. Sato, {\em Toward the classification of 
higher-dimensional toric Fano varieties}, 
Preprint, November 29, 1998. 

\bibitem{S1} Y.-T. Siu, {\em 
The existence of K\"ahler-Einstein metrics on manifolds with positive 
anticanonical line bundle and a suitable finite symmetry group}, 
Ann. Math., II. Ser. {\bf 127}, No.3 (1988), 585-627.  

\bibitem{T1} G. Tian, {\em  
On K\"ahler-Einstein metrics on certain Kaehler manifolds with 
$c_1(M)>0$}, Invent. Math. {\bf 89} (1987), 225-246.  

\bibitem{Tian} G. Tian, {\em 
K\"ahler-Einstein metrics on algebraic manifolds}, 
Proc. Int. Congr. Math., Kyoto/Japan 1990, Vol. I (1991), 587-598.

\bibitem{T3} G. Tian, {\em K\"ahler-Einstein metrics with positive 
scalar curvature}, Invent. Math. 130, No.1 (1997), 1-37. 

\bibitem{TY} G. Tian, S.-T. Yau, 
{\em K\"ahler-Einstein metrics on complex surfaces with $C\sb 1>0$}, 
Commun. Math. Phys. {\bf 112} (1987), 175-203.

\bibitem{VK} V. E. Voskresenskij, A.A.  Klyachko, {\em 
Toroidal Fano varieties and root systems}, 
Math. USSR, Izv. {\bf 24} (1985), 221-244. 


\bibitem{WW} K. Watanabe, M. Watanabe, {\em 
The classification of Fano 3-folds with torus embeddings}, 
Tokyo J. Math. {\bf 5} (1982), 37-48. 

\end{thebibliography}
\end{document}